\documentclass[12pt]{article}

\usepackage{color}
\usepackage{amsmath}
\usepackage{amssymb}
\usepackage{amsthm,url}
\usepackage{dsfont}
\usepackage{textcomp}
\usepackage{wasysym}
\usepackage{graphicx}
\usepackage[shortlabels]{enumitem}

\newtheorem{theorem}{Theorem}
\newtheorem{lemma}[theorem]{Lemma}
\newtheorem{cj}[theorem]{Conjecture}
\newtheorem{proposition}[theorem]{Proposition}

\theoremstyle{definition}
\newtheorem*{remark}{Remark}
\newtheorem*{remark*}{Remark}

\begin{document}

\title{Linking of three triangles in 3-space}

\author{
  E. Kogan
  \footnote{
    We would like to thank A. Skopenkov for helpful discussions and suggestions.
    \newline
    E. Kogan: Higher School of Economics
  }
}
\date{}

\maketitle

\begin{abstract}
Two triples of triangles having pairwise disjoint outlines in 3-space are called {\it combinatorially isotopic} if one triple can be obtained from the other by a continuous motion during which the outlines of the triangles remain pairwise disjoint.
We conjecture that it can be algorithmically checked if an (ordered or unordered) triple of triangles
is combinatorially isotopic to a triple of triangles having pairwise disjoint convex hulls.
We also conjecture that any unordered triple of pairwise disjoint triangles in 3-space
belongs to one of the 5 types of such triples listed in the paper.
We present an elementary proof that triples of different types are not combinatorially isotopic.
\end{abstract}




\section{Introduction and main results}

Examples of linked triangles go back to Middle Ages \cite{Va}.

In this paper, a triangle is the outline (as opposed to the convex hull).

We study the classification of triples of pairwise disjoint triangles up to continuous motion during
which the triangles remain pairwise disjoint.
We conjecture that it can be algorithmically checked if an (ordered or unordered) triple of triangles
is combinatorially isotopic to a triple of triangles having pairwise disjoint convex hulls,
see Conjecture~\ref{unlinked-criterion}.
We also conjecture that any unordered triple of pairwise disjoint triangles in 3-space
belongs to one of the 5 types of such triples listed in Proposition \ref{non-isotopy}, see Conjecture~\ref{main}.
We present an elementary proof that triples of different types are not combinatorially isotopic,
see Proposition~\ref{non-isotopy}.

Triangles are considered as subsets of 3-space, i.e. the vertices are not numbered and the triangles are unoriented.
Moreover, only non-degenerate triangles are considered.

In this paper, when considering an (ordered or unordered) set of triangles,
we require the triangles to be disjoint.
We call such sets \textbf{pairs} and \textbf{triples} of triangles if they contain $2$ and $3$ triangles, respectively.



The {\it convex hull} $\left<X\right>$ of a finite union $X$ of segments in the plane is the smallest convex polygon $P$ in the plane such that $P \supseteq X$.

Let $ABC$ be a triangle from an (ordered or unordered) set $L$ of triangles.
Let $C'$ be a point outside the line $AB$ such that $(\left<ACC'\right>\cup \left<BCC'\right>) \cap \Delta=\emptyset$
for any other triangle $\Delta$ in $L$.
Then an {\bf elementary move} of $L$ is the replacement of $ABC$ by $ABC'$ in $L$.
Two (ordered or unordered) sets of triangles in 3-space are called {\bf combinatorially isotopic} if one can be obtained from the other by a sequence of elementary moves.
This notion was introduced in \cite[\S4.1 `Linking of triangles in space']{Sk} analogously to \cite{PS96} but perhaps was studied earlier.
We conjecture that two (ordered or unordered) sets of triangles in 3-space are combinatorially isotopic iff one set can be obtained from the other by a PL isotopy \cite{isotopy} during which the triangles remain pairwise disjoint triangles.


Obviously, combinatorial isotopy is an equivalence relation.




An (ordered or unordered) pair of triangles is called \textbf{trivial}
if it is combinatorially isotopic to an (ordered or unordered) pair
of triangles whose convex hulls are disjoint.

It is clear that any two trivial (ordered or unordered) pairs of triangles are combinatorially isotopic.
It is also clear that an ordered pair $(\Delta_1, \Delta_2)$ of triangles is trivial
iff the ordered pair $(\Delta_2, \Delta_1)$ is trivial.

An ordered pair $(\Delta_1, \Delta_2)$ of triangles is called \textbf{linked}
(cf. \cite{Sk14})
if the first triangle $\Delta_1$ intersects the convex hull $\left<\Delta_2\right>$ of the second triangle
at exactly one point $P$ and the two segments of $\Delta_1$ going out of $P$
are on the opposite sides of the plane containing $\Delta_2$.
The property of being linked is not symmetric a priori.
However, this is easily proved by considering
$\left<\Delta_1\right> \cap \left<\Delta_2\right>$.

An unordered pair $\{\Delta_1, \Delta_2\}$ of triangles is called \textbf{linked}
if the ordered pairs $(\Delta_1, \Delta_2)$ and $(\Delta_2, \Delta_1)$ are linked.




\begin{theorem}[see proof in \S\ref{proofs}]\label{2-class}
    \leavevmode
    \begin{enumerate}
        \renewcommand{\theenumi}{(\alph{enumi})}
        \renewcommand{\labelenumi}{\theenumi}
            \item\label{2-class-fullness} Any (ordered or unordered) pair of triangles is either trivial or linked.
            \item\label{2-class-non-isotopy} No linked (ordered or unordered) pair of triangles is trivial.
    \end{enumerate}
\end{theorem}


\begin{cj}\label{2-class-nontrivial-correctness}
    Any two linked (ordered or unordered) pairs of triangles are combinatorially isotopic.
\end{cj}

This conjecture is likely to be known.
However, it would be nice if a proof were published.

Consider the following conjecture:
\emph{any ordered pair $(\Delta_1,\Delta_2)$ of triangles
  is combinatorially isotopic to the ordered pair $(\Delta_2,\Delta_1)$.}
This statement follows from Conjecture \ref{2-class-nontrivial-correctness}
and simple Theorem \ref{2-class}\ref{2-class-fullness}.
This shows that Conjecture \ref{2-class-nontrivial-correctness} is not so trivial.




\bigskip

Recall the following examples of unordered triples of triangles.

Take a triple $T^B$ formed by the triangle
$\Delta_0$ with vertices $(2, 0, 0)$, $(-2, \pm 1, 0)$
together with two triangles $\Delta_1$, $\Delta_2$ defined as follows:
$\Delta_j$ is obtained from $\Delta_{j-1}$ by the cyclic permutation $x\to y\to z\to x$ of the coordinates.
An unordered triple of triangles is called \textbf{Borromean} if it is
combinatorially isotopic to the triple $T^B$.

Take a triple $T^3$ formed by an equilateral triangle $\Delta_0$
together with two triangles $\Delta_1$, $\Delta_2$ defined as follows:
$\Delta_j$ is obtained from $\Delta_{j-1}$ by the composition of the rotation
around the altitude $\vec h$ of $\Delta_{0}$ by $\pi/3$, and the translation by $\vec h/3$.
(Cf. \cite[\S4.1 `Linking of triangles in space', 4.1.5]{Sk}.)

\textbf{A triple of type $i$}, $i = 1, 2, 3$
is an unordered triple of triangles such that $i$ pairs of its triangles are linked
and the remaining such pairs are trivial.
It is obvious that there exist triples of types 1 and 2.
The triple $T^3$ is an example of a triple of type 3.


An (ordered or unordered) triple of triangles is called \textbf{trivial} if it is combinatorially isotopic
to an (ordered or unordered) triple of triangles having pairwise disjoint convex hulls.
Clearly if a triple of triangles is trivial,
then all the pairs of the triangles of the triple are trivial.
However, the opposite statement is not correct
by part \ref{bor-neq-trivial} of the following proposition (see also pictures representing counterexamples in \cite{Bo, Va}).
An (ordered or unordered) triple of triangles is called \textbf{3-linked}
if all pairs of the triangles of the triple are trivial but the triple is not trivial.

\begin{proposition}\label{non-isotopy}
    \leavevmode
    \begin{enumerate}[(a)]
            \item\label{bor-neq-trivial} The triple $T^B$ is not trivial.
            \item\label{non-isotopy-corollary} The following types of unordered triples of triangles
        are pairwise distinct:
    trivial, Borromean, type 1, type 2 and type 3.
    \end{enumerate}
\end{proposition}

Proposition \ref{non-isotopy}\ref{non-isotopy-corollary}
is a straightforward corollary of Theorem \ref{2-class} and Proposition \ref{non-isotopy}\ref{bor-neq-trivial}.
We will prove Proposition \ref{non-isotopy}\ref{bor-neq-trivial} after Lemma \ref{bor-isotopy}.

Obviously, an (ordered or unordered) triple of triangles is trivial
iff all the (ordered or unordered) pairs of the triangles of the triple are trivial
and the triple is not 3-linked.
Since the property of an (ordered or unordered) pair of triangles being trivial
can be algorithmically checked,
the following conjecture would give an algorithmic criterion
for the property of a triple of triangles being trivial.

\begin{cj}\label{unlinked-criterion}
    The following properties of an (ordered or unordered) triple $L$ of triangles are equivalent:
    \begin{enumerate}[(i)]
            \item\label{prop:ntriv} $L$ is 3-linked;
            \item\label{prop:3link}
        the triangles of $L$ can be enumerated
        as
        $\Delta_0$, $\Delta_1$, and $\Delta_2$
        so that
        \[
            \left| \Delta_0 \cap \left<\Delta_1\right> \right| =
            \left| \Delta_1 \cap \left<\Delta_2\right> \right| =
            \left| \Delta_2 \cap \left<\Delta_0\right> \right| = 2
        \]
        and $\left<\Delta_0\right> \cap \left<\Delta_1\right> \cap \left<\Delta_2\right>  \neq \emptyset$.
    \end{enumerate}
\end{cj}

The implication \ref{prop:3link} $\Rightarrow$ \ref{prop:ntriv} in Conjecture \ref{unlinked-criterion}
is a straightforward corollary of the following lemma
and the fact
that no (ordered or unordered) triple of triangles with pairwise disjoint convex hulls satisfies \ref{prop:3link}.

\begin{lemma}[see proof in \S\ref{proofs}]\label{bor-isotopy}
    The property \ref{prop:3link} of Conjecture \ref{unlinked-criterion}
    is invariant under combinatorial isotopy.
\end{lemma}

\begin{proof}[Proof of Proposition \ref{non-isotopy}\ref{bor-neq-trivial}]
It is obvious that the triple $T^B$ satisfies the property \ref{prop:3link} of Conjecture \ref{unlinked-criterion}
and that no triple of triangles having disjoint hulls satisfies this property.
Thus the proposition
follows from Lemma \ref{bor-isotopy}.
\end{proof}

In the proof of Lemma \ref{bor-isotopy} in \S\ref{proofs}, we reformulate the property \ref{prop:3link}
of Conjecture \ref{unlinked-criterion}, see Lemma \ref{ii-iii}.




\begin{cj}[Classification Conjecture]\label{main}
    \leavevmode
    \begin{enumerate}[(a)]
            \item\label{main-fullness} Any unordered triple of triangles
        belongs to one of the types of triples listed in Proposition \ref{non-isotopy}\ref{non-isotopy-corollary}.
            \item\label{main-correctness} Triples of the same type are combinatorially isotopic.
    \end{enumerate}

\end{cj}

This conjecture should be compared to Milnor's {\it link homotopy} classification of \textit{3-component links} \cite{Mi54},
see statement accessible to non-specialists in \cite[\S 4.4 `Borromean rings and commutators', Milnor's theorem 4.4.3]{Sk}.

The versions of Proposition \ref{non-isotopy}\ref{non-isotopy-corollary}
and Conjecture \ref{main}
for ordered triples of triangles are analogous
with 6 types of ordered triples instead of types 1 and 2.

The following conjecture follows from Conjecture \ref{main}\ref{main-fullness}:
{\it any unordered 3-linked triple of triangles is combinatorially isotopic to $T^B$.}



\begin{remark}[not used in the sequel]
    A triple of triangles satisfies the property \ref{prop:3link} of Conjecture \ref{unlinked-criterion}
    iff its Milnor-Massey number is 1.
    Thus Lemma \ref{bor-isotopy}
    can be obtained using the Milnor-Massey number.
    See an elementary definition of \emph{the Milnor-Massey number} in
    \cite[\S 4.6 `Triple linking modulo~2', definition after problem 4.6.5]{Sk}.
\end{remark}





\section{Proofs}\label{proofs}

\begin{proof}[Proof of Theorem \ref{2-class}\ref{2-class-fullness}]
        The case of unordered pairs of triangles follows from the case of ordered pairs of triangles
        because if an ordered pair $(\Delta_1, \Delta_2)$
        is trivial (linked), then the unordered pair $\{\Delta_1, \Delta_2\}$ is trivial (linked).
        So let us prove Theorem \ref{2-class}\ref{2-class-fullness} for the case of ordered pairs of triangles.

        Denote by $L = (\Delta_1, \Delta_2)$ an ordered pair of triangles.
        Suppose that $L$ is not linked. We need to prove that $L$ is trivial.

        Firstly, either the intersection $\left<\Delta_1\right> \cap \Delta_2$
        consists of two points or $\Delta_2$ lies in one half-space with respect to the plane containing $\Delta_1$.
        Thus we can perform elementary moves by a small distance so
        that the intersection $\left<\Delta_1\right> \cap \Delta_2$
        either stays the same if it consists of two points
        or vanishes otherwise.
        Denote the obtained pair of triangles by $L' = (\Delta'_1, \Delta'_2)$.
        It is sufficient to prove that $L'$ is trivial.


        \emph{The case where $\left<\Delta'_1\right> \cap \Delta'_2 = \emptyset$.}
        There is a sequence of elementary moves of $L'$
        yielding an ordered pair $(\Delta''_1, \Delta'_2)$ of triangles
        such that the triangle $\Delta''_1$
        lies in the same plane as $\Delta'_1$ and
        does not intersect the set
        $\left<\Delta'_1\right> \cap \left<\Delta'_2\right>$.
        The convex hulls
        $\left<\Delta''_1\right>$ and $\left<\Delta'_2\right>$ do not intersect.
        Thus the resulting pair $(\Delta''_1, \Delta'_2)$ is trivial.
        Hence the pair $L'$ is trivial as well.


        \emph{The case where $\left|\left<\Delta'_1\right> \cap \Delta'_2\right| = 2$.}
        In this case
        the outline of $\Delta'_1$ does not intersect the convex hull of $\Delta'_2$
        and hence the pair $(\Delta'_2, \Delta'_1)$ is trivial by the first case.
        Then the pair $L'$ is trivial as well.

\end{proof}

\begin{proof}[Proof of Theorem \ref{2-class}\ref{2-class-non-isotopy}]
    It is sufficient to prove that the property of being linked
    is invariant under combinatorial isotopy.
    We will prove this in the case of an ordered pair. The unordered case is analogous.

    Consider an elementary move $(\Delta, ABC) \to (\Delta, ABC')$ of a pair of triangles.
    The case when the first triangle is changed is analogous by the symmetry of the property of being linked.

    An oriented triangle $\Delta$ `goes inside' the tetrahedron $ABCC'$
    as many times as it `goes outside' the tetrahedron.
    Then this triangle `crosses' the face $\left<ABC'\right>$
    at an odd number of points which is less than $2$.
    Thus the pair $(\Delta, ABC')$ is linked.
\end{proof}

\begin{lemma}[see proof in \S\ref{proofs}]\label{ii-iii}
    Let $L$ be an (ordered or unordered) triple of triangles.
    The property \ref{prop:3link} of Conjecture \ref{unlinked-criterion}
    is equivalent to the following:

    \begin{enumerate}[({ii}i)]
            \item\label{prop:3'link} the triangles of $L$ can be enumerated
        as
        $\Delta_0$, $\Delta_1$, and $\Delta_2$
        so that
        \[
            \Delta_1 \cap \left<\Delta_0\right> = \emptyset,\quad
            \left| \Delta_1 \cap \left<\Delta_2\right> \right| =
            \left| \Delta_2 \cap \left<\Delta_0\right> \right| = 2
        \]
        and $\left<\Delta_0\right> \cap \left<\Delta_1\right> \cap \left<\Delta_2\right>  \neq \emptyset$.
    \end{enumerate}
\end{lemma}

\begin{proof}[Proof of Lemma \ref{ii-iii}]

The implication \ref{prop:3link} $\Rightarrow$ \ref{prop:3'link} is obvious.
Let us prove the implication \ref{prop:3'link} $\Rightarrow$ \ref{prop:3link}.
To do this, we need to prove that $\left| \Delta_0 \cap \left<\Delta_1\right> \right| = 2$.

\begin{figure}[h]
\centering
\includegraphics[scale=0.2]{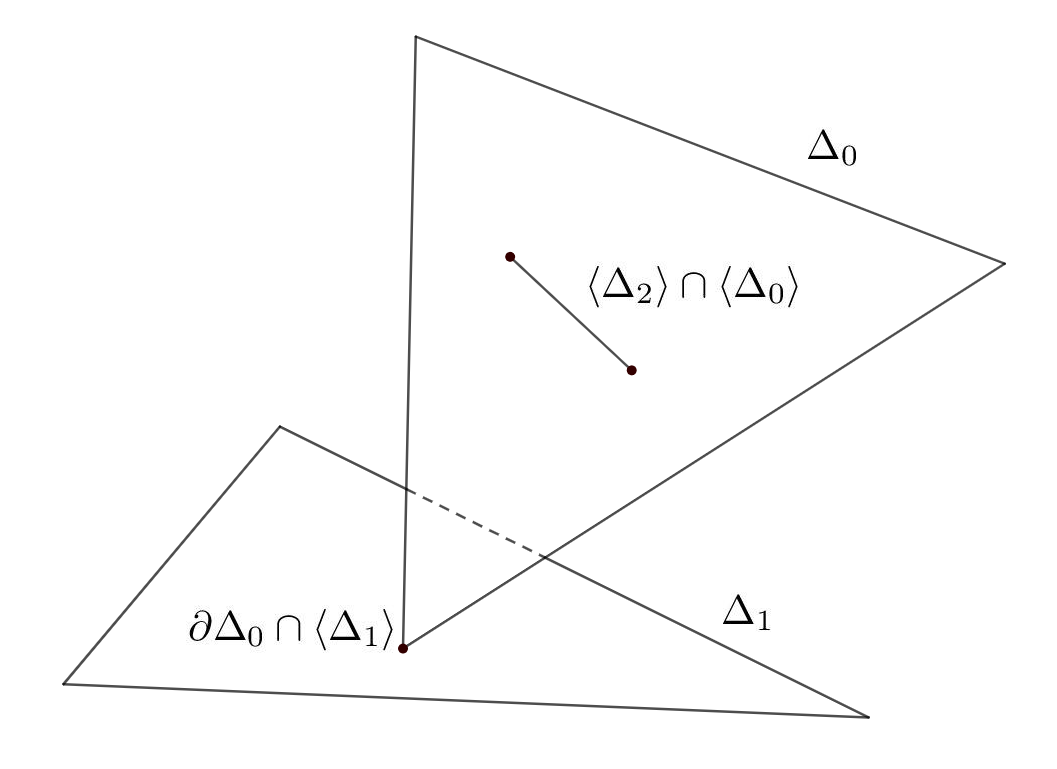}
\caption{
	$\left| \Delta_0 \cap \left<\Delta_1\right> \right| = 1$
}
\label{L0}
\end{figure}


We have $\left<\Delta_0\right> \cap \left<\Delta_1\right> \neq \emptyset$.
Hence
\begin{equation*}
    (\Delta_0 \cap \left<\Delta_1\right>) \cup (\Delta_1 \cap \left<\Delta_0\right>)
    \supset
    \partial(\left<\Delta_0\right> \cap \left<\Delta_1\right>) \neq \emptyset
\end{equation*}
Since $\Delta_1 \cap \left<\Delta_0\right> = \emptyset$,
it follows that $\Delta_0 \cap \left<\Delta_1\right> \neq \emptyset$.

Therefore $\Delta_0$ intersects $\left<\Delta_1\right>$
at a point, or at two points,
or at a segment, or at the whole triangle $\Delta_0$.

\emph{The case when $\Delta_0 \cap \left<\Delta_1\right> = \Delta_0$.}
In this case, $\Delta_0 \subset \left<\Delta_1\right>$.
Since $\left| \Delta_1 \cap \left<\Delta_2\right> \right| = 2$,
we have $\Delta_2 \cap \left<\Delta_1\right> = \emptyset$.
This and $\left<\Delta_0\right> \subset \left<\Delta_1\right>$
imply that
$$
\emptyset
= \Delta_2 \cap \left<\Delta_1\right> \supset
\Delta_2 \cap \left<\Delta_0\right>
\cap \left<\Delta_1\right>
= \Delta_2 \cap \left<\Delta_0\right> \neq \emptyset.
$$
A contradiction.

\emph{The case when $\Delta_0 \cap \left<\Delta_1\right>$ is a point or a segment (see Fig. \ref{L0}).}
In this case, $\Delta_0$ is in the same half-space with respect to the plane $l(\Delta_1)$ containing $\Delta_1$.
Then the intersection $\left<\Delta_2\right> \cap \left<\Delta_0\right>$
is the segment between the two points
of $\Delta_2 \cap \left<\Delta_0\right>$.
These two points are in the same half-space with respect to the plane $l(\Delta_1)$
and outside of the plane.
Hence $\left<\Delta_2\right> \cap \left<\Delta_0\right>$ does not intersect $l(\Delta_1)$.
This is a contradiction
because $\left<\Delta_2\right> \cap \left<\Delta_0\right>$
intersects $\left<\Delta_1\right> \subset l(\Delta_1)$.

Thus none of these two cases can take place and $\Delta_0$ intersects $\left<\Delta_1\right>$ at two points.

\end{proof}

\begin{proof}[Proof of Lemma \ref{bor-isotopy}]

Suppose that $T$ is a 3-linked triple of triangles $\Delta_0, \Delta_1, \Delta_2$
such that $\left| \Delta_j \cap \left<\Delta_{(j + 1)\!\!\mod{\!3}}\right> \right| = 2$.
Suppose further that we are given an elementary move of $T$
in which the triangle $\Delta_0 = ABC$ is replaced
by the triangle $\Delta_0' = ABC'$.
It is sufficient to prove that
the triple $\left(\Delta_0', \Delta_1, \Delta_2\right)$
satisfies the property \ref{prop:3'link} of Lemma \ref{ii-iii}.

Denote by $P$ and $Q$
the two points of the intersection $\Delta_1 \cap \left<\Delta_2\right>$.
Then $\left<\Delta_1\right> \cap \left<\Delta_2\right>$
is the segment $PQ$.
Since $\left<\Delta_0\right> \cap \left<\Delta_1\right> \cap \left<\Delta_2\right> \neq \emptyset$, we have
$PQ \cap \left<\Delta_0\right> \neq \emptyset$.
Hence $P$ and $Q$ are on different sides of the plane containing $\Delta_0$.
Without loss of generality $P$ and $C'$ are on the same side of this plane.

Denote by $\tau$ the tetrahedron $ABCC'$.

\begin{figure}[h]
\centering
\includegraphics[scale=0.2]{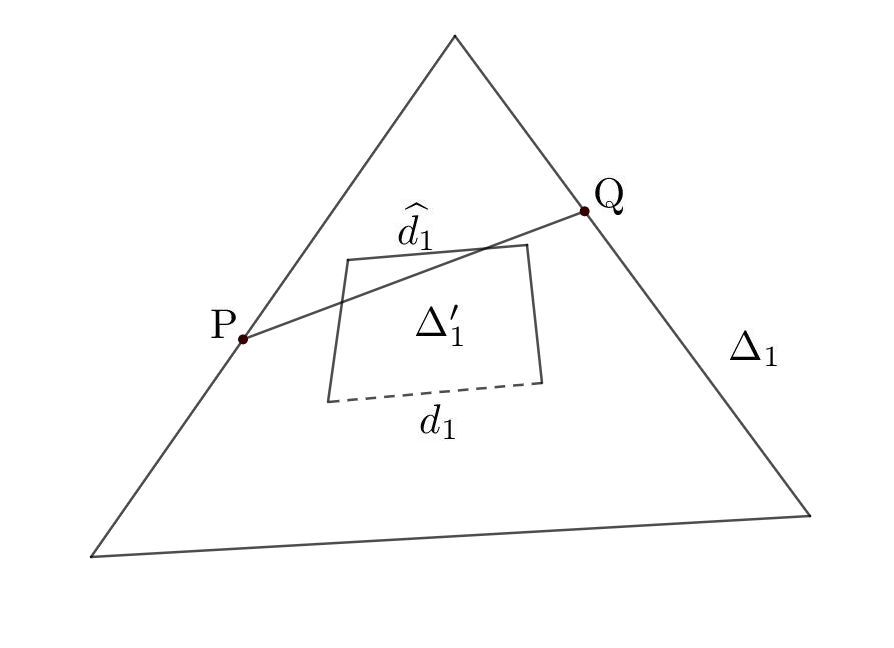}
\caption{$\Delta_1$ and $\Delta_1'$}
\label{D1}
\end{figure}


In the following paragraph we prove that $\Delta_1 \cap \left<\Delta_0'\right> = \emptyset$ and that $P$ and $Q$ are in the exterior of $\tau$.

Let $\Delta_1'$ be the intersection of $\tau$ and the plane containing $\Delta_1$ (see Fig.~\ref{D1}).
Out of the faces of $\tau$, the triangle $\Delta_1$ may only intersect $ABC'$.
Hence $\Delta_1$ cannot intersect more than one side of $\Delta_1'$.
Denote by $d_1$ the side of $\Delta_1'$ that $\Delta_1$ intersects
or any side of $\Delta_1'$ not contained in $\left<\Delta_0\right>$
if $\Delta_1$ does not intersect $\Delta_1'$.
Denote $\widehat{d_1} := \Delta_1' - d_1$.
The segment $PQ$ intersects $\widehat{d_1} \cap \left<\Delta_0\right>$,
hence $\emptyset \neq PQ \cap \widehat{d_1} \subset \left <\Delta_1 \right >$.
Consequently, $\widehat{d_1}$
is in the interior of the triangle $\left<\Delta_1\right>$.
Hence the polygon $\Delta_1'$ is in the interior of $\left<\Delta_1\right>$,
so $\Delta_1$ does not intersect $\Delta_1'$.
Therefore $\Delta_1$ does not intersect either $\tau$
or $\left<\Delta_0'\right>$.
Also, since $\Delta_1'$ is in the interior of $\left<\Delta_1\right>$,
the points $P$ and $Q$ are in the exterior of $\tau$.


In the following paragraph we prove that
$\left| \Delta_2 \cap \left<\Delta_0'\right> \right| = 2$.

Denote by $\Delta_2^P$ the part of $\Delta_2$ which is
on the same side of the plane containing $\Delta_0$ as $C'$ (and $P$).
Then $P$ is in the interior of $\Delta_2^P$ and in the exterior of $\tau$.
However, in all cases where
the intersection $\Delta_2 \cap \left<\Delta_0\right>$
does not consist of two points
(whether it is an empty set, or a point, or a segment),
$\Delta_2^P$ is entirely in the interior of $\tau$
which contradicts the existence of $P$.
Thus $\left| \Delta_2 \cap \left<\Delta_0'\right> \right| = 2$.

Finally we show that $\left<\Delta_0'\right> \cap \left<\Delta_1\right> \cap \left<\Delta_2\right> \neq \emptyset$.

Denote by $\Delta_2'$ the intersection of $\tau$
and the plane containing $\Delta_2$.
Then $\Delta_2$ intersects two sides of $\Delta_2'$ at two points each.
Hence the other side(s) of $\Delta_2'$ is (are) in the exterior of $\left<\Delta_2\right>$.
Then the faces $\left<ACC'\right>$ and $\left<BCC'\right>$ of the tetrahedron $\tau$
do not intersect $\left<\Delta_2\right> \supset PQ$.
Moreover, since $P$ and $Q$ are in the exterior of $\tau$
and $PQ \cap \tau \supset PQ \cap \left<ABC\right> \neq \emptyset$,
the segment $PQ$ intersects two faces of $\tau$.
Hence $PQ$ intersects $\left<ABC\right> = \left<\Delta_0\right>$
and $\left<ABC'\right> = \left<\Delta_0'\right>$,
which means that $\left<\Delta_1\right> \cap \left<\Delta_2\right> \cap \left<\Delta_0'\right> \neq \emptyset$.
\end{proof}

\end{document}